\documentclass[11pt,a4]{article}

\usepackage[mathscr]{eucal}
\usepackage{amsmath}
\usepackage{graphicx}
\usepackage{amssymb}
\usepackage{latexsym}
\usepackage{amsthm}
\usepackage{amscd}
\usepackage{enumerate}
\theoremstyle{plain}
\newtheorem{thm}{Theorem}[section]

\theoremstyle{definition}
\newtheorem{Def}[thm]{Definition}

\title{Generalizations of Sylvester's refinement of Euler's odd-strict theorem
}
\author{MASANORI ANDO\ \ (Nara Gakuen University)}
\date{}
\begin{document}
\pagestyle{empty}
\maketitle\thispagestyle{empty}
\section{Introduction}
Regarding Euler's odd-strict theorem, which is the most basic partition identity,
A refinement was done by Sylvester\cite{S}, and it was generalized by Bessenrodt to the $r$-regular and $r$-class regular cases\cite{B}.
In this paper, we focus on the periodicity of the exponent seen on the $r$-regular side in Bessenrodt's generalization, and further generalize from that point of view. 
At the end of this paper, we also introduce a partition identity when extended without regard to the periodicity of the exponent.
%最も基本的な分割恒等式である, Euler のodd-strict 定理について, 
%Sylvester による精密化が行われ\cite{S}, さらにそれがBessenrodt によって$r$-正則, $r$-類正則の場合に一般化されたcite{B}. 
%本論文においては, Bessenrodt による一般化の際に$r$-正則側に見られる指数の周期性に着目し, その観点からの一般化を行う. 
%また最後に, 証明に用いた写像から対応が付く範囲いっぱいまで拡張した場合の分割恒等式も紹介する. 
\section{Notation}
\begin{Def}
Let $n$ be a positive integer. A partition $\lambda $ of $n$ is an integer sequence 
\[
\lambda =(\lambda_1,\lambda_2,\ldots,\lambda_\ell) 
\]
satisfying $\lambda _1 \geq \lambda _2 \geq \ldots \geq \lambda _\ell >0$ and $\displaystyle \sum _{i=1}^{\ell }{\lambda _i} =n$. 
We call $\ell (\lambda):= \ell$ the length of $\lambda $, 
$|\lambda|:=n$ the size of $\lambda $, 
and each $\lambda _i$ a part of $\lambda $. 
We let ${\mathcal {P}} $ and ${\mathcal {P}}(n) $ denote the set of partitions and the set of partitions of $n$. 
From now on, $(n)$ will represent ``size $n$''. 
We denote the number of part $k$ in $\lambda$ as $m_k(\lambda)$. 
We also represent $\lambda =(\cdots 2^{m_2(\lambda)}1^{m_1(\lambda)})$. 
\end{Def}
\noindent
\textbf{Example.} $n=5$
\begin{eqnarray*}
\mathcal{P}(5)&=&\{
(5), (4,1), (3,2), (3,1,1), (2,2,1), (2,1,1,1), (1,1,1,1,1)
\}\\
&=&\{
(5), (41), (32), (31^2), (2^21), (21^3), (1^5)
\}
\end{eqnarray*}
\begin{Def}
Let $r$ be a positive integer. 
We define two subsets of $\mathcal{P}$ as follows. 
\[
\mathcal{RP}_r:=
\{
\lambda \in \mathcal{P}\ |\ {}^\forall k, m_k(\lambda)<r
\}
\]
\[
\mathcal{CP}_r:=
\{
\lambda\in \mathcal{P}\ |\ {}^\forall k, m_{rk}(\lambda)= 0
\}
\]
We call them $r$-regular partitions and $r$-class regular partitions. 
\end{Def}
\noindent
\textbf{Example.} $n=6, r=3$
\[
\mathcal{RP}_3(6)=\{(6), (51),(42),(41^2), (3^2), (321), (2^21^2)\}
\]
\[
\mathcal{CP}_3(6)=\{(51), (42), (41^2), (2^3), (2^21^2), (21^4), (1^6)\}
\]
It is well known that $\sharp\mathcal{RP}_r(n)=\sharp\mathcal{CP}_r(n)$. 
\begin{Def} 
%成分を$r$で割った商と余りの記号を定める. 
We denote the quotient and remainder symbols when the part is divided by $r$.
\[
Q_r(\lambda_i):=\left\lfloor\frac{\lambda_i}{r}\right\rfloor, 
R_r(\lambda_i):=\lambda_i-rQ_r(\lambda_i), 
QR_r(\lambda_i):=(Q_r(\lambda_i), R_r(\lambda_i))
\]
\end{Def}
\section{refined version}
%成分の種数と連続整数列の個数
We denote the number of kind of parts in $\lambda$ as $K(\lambda)$ and the number of consecutive integer sequences in $\lambda$ as $S(\lambda)$. 
\[
K(\lambda):=\sharp\{k\ |\ m_{k}(\lambda)>0\}
, 
S(\lambda):=\sharp\{i\ |\ \lambda_i>\lambda_{i+1}+1\}+1
\]
\begin{thm}[Sylvester\cite{S}]
%ちょうど$m$種の奇数成分による$n$の分割の個数は, 相異なる成分への$n$の分割で, 独立した$m$個の連続正整数列からなるものの個数に等しい. 
The number of odd partitions of $n$ by just $m$ kind of odd parts is equal to 
the number of strict partitions of $n$ by independent $ m $ consecutive positive integer sequences.
\[
\sharp\left\{
\lambda\in \mathcal{CP}_2(n)\ |\ K(\lambda)=m
\right\}
=
\sharp\left\{
\lambda\in \mathcal{RP}_2(n)\ |\  S(\lambda)=m
\right\}
\]
\end{thm}
%Sylvester のこの定理はBessenrodt によって次の形に一般化されている. 
This Sylvester's theorem is generalized by Bessenrodt to the following form.
\begin{thm}[Bessenrodt\cite{B}]
For any positive integer $r, 1\leq j\leq r-1$, 
The number of partitions of $n$ such that all parts are congruent with $j$ modulo $r$ is equal to
 the number of partitions of $n$ such that $\lambda =(k_1^jk_2^{r-j}k_3^jk_4^{r-j}\ldots), k_1>k_2>k_3>k_4>\ldots$. 
\begin{eqnarray*}
\mathcal{CP}_{\{r,j\}}(n)&:=&
\{
\lambda \in \mathcal{P}(n)\ |\ {}^\forall i\leq \ell(\lambda), R_r(\lambda_i)=j
\}\\
\mathcal{RP}_{\{r,j\}}(n)&:=&
\{
(k_1^jk_2^{r-j}k_3^jk_4^{r-j}\ldots) \in \mathcal{P}(n)\ |\ k_1>k_2>k_3>k_4>\ldots> 0
\}
\end{eqnarray*}
\[
\sharp\mathcal{CP}_{\{r,j\}}(n)=\sharp\mathcal{RP}_{\{r,j\}}(n)
\]
%任意の自然数$r$, $1\leq j\leq r-1$ について, 
%$ri+j$型成分による$n$の分割の個数は, 
%大きさ$n$の$r$-正則分割で, $\lambda =(k_1^jk_2^{r-j}k_3^jk_4^{r-j}\ldots), k_1>k_2>k_3>k_4>\ldots$,  
%の形のものの個数に等しい. ]
Especially, for all $m$, 
\[
\sharp\{\lambda\in\mathcal{CP}_{\{r,j\}}(n)\ |\ K(\lambda)=m\}=\sharp\{\lambda\in\mathcal{RP}_{\{r,j\}}(n)\ |\ S(\lambda)=m\}
\]
\end{thm}
\noindent
\textbf{Example1. } Let $n=8, r=3, j=1$. The set of partitions of $8$ that all parts congruent to $1$ modulo $3$ is as follows.
%\textbf{例.} $n=8, r=3, j=1$ の場合, $3$を法として$1$ に合同な成分のみによる$8$の分割は, 
\[
\mathcal{CP}_{\{3,1\}}(8)=
\{
(71), (4^2), (41^4), (1^8)
\}
\]
And the set of partitions on $8$ with exponent 121212... is next.
%一方指数が$1212\ldots $ となるような$8$ の分割は, 
\[
\mathcal{RP}_{\{3,1\}}(8)=
\{
(8^1), (6^11^2), (4^12^2), (3^12^21^1)
\}
\]
The case of $j=2$, the set of partitions of $8$ that all parts congruent to $2$ modulo $3$ is 
%$j=2$の場合, $3$を法として$2$ に合同な成分のみによる$8$の分割は, 
\[
\mathcal{OP}_{3,2}(8)=
\{
(8), (2^4)
\}
\]
The set of partitions on $8$ with exponent 212121... is 
%指数が$2121\ldots $ となるような$8$ の分割は, 
\[
\mathcal{SP}_{3,2}(8)=
\{
(4^2), (3^22^1)
\}
\]
In both cases the order of the set is equal.
The equality also holds under the conditions $K(\lambda), S(\lambda)$. 
%成分の種数と連続整数列の個数との関係についても正しいことが確認できる. 
\section{Repeating Regular Partition}
%Bessenrodt の定理は分割の指数の周期が2の場合であるため, 
In Bessenrodt's theorem, it can be seen that the exponent of the partitions of $\mathcal{RP}$ circulates with a period of 2.
%これを一般に周期的な場合へと拡張する. 
We extend this to the general periodic case.

%周期2の場合であれば, 成分0も許すと思えば, 分割の長さはrの倍数の場合のみを考えていることになる. 
In the theorem for the period 2 case, If we think that part 0 is also allowed, the length of the partition is a multiple of $2$. 
Three generalizations can be considered: the case where the length of the partition is free, the case where the length of the partiton is a multiple of $r$, and the case where the length of the partition is a multiple of $r$ with the component 0 allowed.

First, the case where the length of the partition is free. 
%任意の$r, t, $
%$1\leq s_1<s_2<\ldots <s_{t-1}< r$ について, 
%大きい成分から順にその個数が$s_1, s_2-s_1, \ldots , r-s_{t-1}$ と周期的になっている分割を考えよう. 
%すなわち, 
For all $r, t$, and $1\leq s_1<s_2<\ldots <s_{t-1}<r$, 
Consider the case that the index of the partiton is repeating $s_1, s_2-s_1, \ldots , r-s_{t-1}$ in order from the largest part. 
%That is, 
%\[
%{}^\forall k, 0\leq{}^\forall j\leq t, m_{\lambda_{rk+s_j}}(\lambda )=s_j-s_{j-1}, s_0=0, s_t=r. 
%\]
%そのような分割全体のなす集合を
%$\mathcal{RP}_{\{r; s_1, s_2, \ldots , s_{t-1}\}}$ と表すこととする. 
We define such a set of partitions $\mathcal{RP}_{\{r; s_1, s_2, \ldots , s_{t-1}\}}$. 
\begin{eqnarray*}
&&\mathcal{RP}_{\{r; s_1, s_2, \ldots , s_{t-1}\}}\\
&:=&
\left\{
\lambda\in \mathcal{RP}_r\ |\ 
{}^\forall k, 0<{}^\forall j\leq t, m_{\lambda_{rk+s_j}}(\lambda )=s_j-s_{j-1}, s_0=0, s_t=r. 
\right\}\\
&=&\left\{
\lambda=({k_1}^{s_1}{k_2}^{s_2-s_1}\ldots {k_t}^{r-s_{t-1}}{k_{t+1}}^{s_1}{k_{t+2}}^{s_2-s_1}\ldots)\ |\ k_1>k_2>\ldots
\right\}
\end{eqnarray*}
%また対応する集合として次を定義する. 
And we define corresponding next set. 
\begin{eqnarray*}
&&\mathcal{CP}_{\{r;s_1, s_2, \ldots , s_{t-1}\}}\\
&:=&
\left\{
\lambda\in \mathcal{CP}_r
\ {\Bigg{|}}\ 
\begin{array}{ll}
{}^\forall i, 0<{}^\exists j<t, R_r(\lambda_i)= s_j, \\
Q_r(\lambda_i)\geq i \Rightarrow R_r(\lambda_i)= s_{t-1}, 
{}^\forall a<t-1, {}^\exists k, 
QR_r(\lambda_k)=(i-1, s_a)\\
%\lambda_k=r(i-1)+s_a \\
Q_r(\lambda_i)=i-1\wedge R_r(\lambda_i)=s_j\Rightarrow {}^\forall a<j, {}^\exists k, 
QR_r(\lambda_k)=(i-1, s_a)
\end{array}
\right\}
\end{eqnarray*}

%完全循環ver. 
%実は完全循環バージョンはSylvester, Bessenrodt の一般化ではない. $r=2$や周期２に特殊化した場合にはそれらの精密化となっている. 
The case where the length of the partiton is a multiple of $r$ and the case where the length of the partition is a multiple of $r$ with the component 0 allowed, 
The sets are subset of First case. 
\[
\mathcal{RP}'_{\{r; s_1, s_2, \ldots , s_{t-1}\}}:=
\left\{
\lambda\in \mathcal{RP}_{\{r; s_1, s_2, \ldots , s_{t-1}\}}\ |\ 
\ell(\lambda)\equiv 0({\rm{mod}}\ r), 
\right\}
\]
\[
\mathcal{CP}'_{\{r;s_1, s_2, \ldots , s_{t-1}\}}:=
\left\{
\lambda\in \mathcal{CP}_{\{r; s_1, s_2, \ldots , s_{t-1}\}}\ |\ 
\begin{array}{ll}
{}^\forall i, Q_r(\lambda_i)\not=i-1\\
\end{array}
\right\}
\]
\[
\mathcal{RP}''_{\{r; s_1, s_2, \ldots , s_{t-1}\}}:=
\left\{
\lambda\in \mathcal{RP}_{\{r; s_1, s_2, \ldots , s_{t-1}\}}\ |\ 
\ell(\lambda)\equiv 0, s_{t-1}({\rm{mod}}\ r), 
\right\}
\]
\[
\mathcal{CP}''_{\{r;s_1, s_2, \ldots , s_{t-1}\}}:=
\left\{
\lambda\in \mathcal{CP}_{\{r; s_1, s_2, \ldots , s_{t-1}\}}\ |\ 
\begin{array}{ll}
{}^\forall i, 
Q_r(\lambda_i)\geq i-1 \Rightarrow R_r(\lambda_i)= s_{t-1}
\end{array}
\right\}
\]
%ここで$I_r(\lambda)=\max\{i \ |\ \lambda_i>r(i-1)\}$
\begin{thm}
For all $n, r, t, 1\leq s_1<s_2<\ldots <s_{t-1}<r$, 
%任意の$n, r, t, $
%$1\leq s_1<s_2<\ldots <s_{t-1}<r$ について, 
\begin{eqnarray*}
\sharp\mathcal{CP}_{\{r;s_1, s_2, \ldots , s_{t-1}\}}(n)=
\sharp\mathcal{RP}_{\{r;s_1, s_2, \ldots , s_{t-1}\}}(n), \\
\sharp\mathcal{CP}'_{\{r;s_1, s_2, \ldots , s_{t-1}\}}(n)=
\sharp\mathcal{RP}'_{\{r;s_1, s_2, \ldots , s_{t-1}\}}(n), \\
\sharp\mathcal{CP}''_{\{r;s_1, s_2, \ldots , s_{t-1}\}}(n)=
\sharp\mathcal{RP}''_{\{r;s_1, s_2, \ldots , s_{t-1}\}}(n)\\
\end{eqnarray*}
\end{thm}
\noindent
\textbf{Example2. }\ For $n=10, r=3, t=3, (s_1, s_2)=(1, 2)$,  
\[
\mathcal{RP}_{\{3;1,2\}}(10)=\left\{ 
\begin{array}{l}
(10^1), (9^11^1), (8^12^1), (7^13^1), (7^12^11^1), (6^14^1), \\
(6^13^11^1), (5^14^11^1), (5^13^12^1), (4^13^12^11^1)
\end{array}
\right\}
\]
\[
\mathcal{CP}_{\{3;1,2\}}(10)=\{ 
(81^2), (541), (52^21), (521^3), (51^5), (2^41^2), (2^31^4), (2^21^6), (21^8), (1^{10})
\}
\]
Then, $\sharp\mathcal{SP}_{\{3;1,2\}}(10)=
\sharp\mathcal{OP}_{\{3;1,2\}}(10)=10$. \\
\textbf{Example2. }\ For $n=12, r=6, t=3, (s_1, s_2)=(1, 3)$, 
\[
\mathcal{RP}_{\{6;1,3\}}(12)=\{ 
(12^1), (10^11^2), (8^12^2), (6^13^2), (5^12^21^3)
\}
\]
\[
\mathcal{CP}_{\{6;1,3\}}(12)=\{
(91^3), (3^31^3), (3^21^6), (31^9), (1^{12})
\}
\]
\textbf{Remark. }
For example, $(93), (731^2), (3^4)\not\in\mathcal{CP}_{\{6;1,3\}}(12)$\\
%\textbf{例4. }\ $n=13, r=6, t=3, (s_1, s_2)=(1, 3)$ について, 
%\[
%\mathcal{RP}_{\{6;1,3\}}(13)=\{ 
%(13^1), (11^11^2), (9^12^2), (7^13^2), (5^14^2)(6^12^21^3), (4^13^21^3)
%\}
%\]
%\[
%\mathcal{CP}_{\{6;1,3\}}(13)=\{
%(931), (91^4), (3^41), (3^31^4), (3^21^7), (31^{10}), (1^{13})
%\}
%\]
\begin{proof}
First we present the Bessenrodt's map. 
It is the same as for period 2 and general period. 
And we construct it's inverse map. 
%先の定理の証明となるBessenrodt の対応を紹介する. 
%同様の対応を一般の周期性の場合に広げたものについて, 逆写像を構成する. 
%構成の仕方から逆写像であることは明らかかと思う. 
For $\mathcal{CP}_{\{r;s_1,s_2, \ldots , s_{t-1}\}}(n)\ni \lambda=\{ \lambda_1 \lambda_2 \ldots \lambda_\ell\}$, 
we correspond the tableau $T(\lambda)$. \\
%$\mathcal{CP}_{\{r;s_1,s_2, \ldots , s_{t-1}\}}(n)\ni \lambda=\{ \lambda_1 \lambda_2 \ldots \lambda_\ell\}$
%について, 次の盤$T(\lambda)$(仮称)を対応させる. \\
{\textbf{Example5. }}\ For $\displaystyle \mathcal{CP}_{\{ 6;1,3\}}(35)\ni\lambda=(15,9,7,3,1)$, 
\[
T(\lambda)=
\begin{array}{ccc}
6&6&3\\
6&3&\\
6&1&\\
3&&\\
1&&\\
\end{array}
\]
When $\lambda_k=ri+s_j$, we put $i$ $r$s and one $s_j$ on $k$-th row. 
%すなわち, $\lambda_k=ri+s_j$ であるとき, $k$行目には, $i$個の$r$ と1個の$s_j$ を置く. 
Becouse of definition of partition, the shape of tableau is Young diagram. 
And each number is weakly decreasing from the left and top. 
For this tableau, we constract $S(\lambda)\in\mathcal{RP}_{\{r;s_1,s_2, \ldots , s_{t-1}\}}(n)$ 
by removing the diagonal hook from $T(\lambda)$ as follows.\\
%このとき分割の定義から, 盤の形はヤング図形であり, 各数字は, 左から右に広義単調減少（右端のみ減少）, 
%上から下に広義単調減少となっている. 
%この盤に対して, 順に対角フックを抜く操作によって, $\mathcal{RP}$の分割を次の様に対応させる。\\
{\textbf{Example6. }}\ 
\[
\begin{array}{ccc}
6&6&3\\
6&3&\\
6&1&\\
3&&\\
1&&\\
\end{array}
\mapsto_7
\begin{array}{ccc}
5&5&2\\
5&3&\\
5&1&\\
2&&\\
&&\\
\end{array}
\mapsto_6
\begin{array}{ccc}
4&4&1\\
4&3&\\
4&1&\\
1&&\\
&&\\
\end{array}
\mapsto_6
\begin{array}{ccc}
3&3&\\
3&3&\\
3&1&\\
&&\\
&&\\
\end{array}
\mapsto_4
\begin{array}{ccc}
2&2&\\
2&3&\\
2&1&\\
&&\\
&&\\
\end{array}
\]
\[
\mapsto_4
\begin{array}{ccc}
1&1&\\
1&3&\\
1&1&\\
\end{array}
\mapsto_4
\begin{array}{cc}
3&\\
1&\\
\\
\end{array}
\mapsto_2
\begin{array}{cc}
2&\\
&\\
\\
\end{array}
\mapsto_1
\begin{array}{cc}
1&\\
&\\
\\
\end{array}
\mapsto_1
\emptyset
\]
Then, $S(\lambda)=(7^16^24^32^11^2)$. 
The exponent of $S(\lambda)$ is repeats $123$, so $S(\lambda)\in \mathcal{RP}_{\{6;1,3\}}(35)$. 
%指数を見ると, $123$ が循環しているため, これは$\mathcal{RP}$の分割である. 

Next we constract inverse map of $S$. \\
 For $\mu=(l_1^{s_1}l_2^{s_2-s_1}\ldots l_t^{r-s_{t-1}}l_{t+1}^{s_1}l_{t+2}^{s_2-s_1}\ldots l_{2t}^{r-s_{t-1}}\ldots)\in \mathcal{RP}_{\{r;s_1,s_2, \ldots , s_{t-1}\}}(n)$, 
we separate its parts by $r$. \\
Let $\lambda^1=(l_1^{s_1}l_2^{s_2-s_1}\ldots l_t^{r-s_{t-1}})$, $\lambda^2=(l_{t+1}^{s_1}l_{t+2}^{s_2-s_1}\ldots l_{2t}^{r-s_{t-1}})$, ....
For each $\lambda^i$, we consider the conjugate ${}^t(\lambda^i)$
\footnote{It's confusing, but the upper left $t$ means transpose.}. 
The parts of ${}^t(\lambda^i)$ belong to $\{s_1, s_2, \ldots , s_{t-1}\}$, $s_t=r$. 
And $m_{s_j}({}^t({\lambda}^i))\geq 1$. 
%一方逆写像について, $\mathcal{RP}\ni\mu=(l_1^{s_1}l_2^{s_2-s_1}\ldots l_t^{r-s_{t-1}}l_{t+1}^{s_1}l_{t+2}^{s_2-s_1}\ldots l_{2t}^{r-s_{t-1}}\ldots)$ に対して, 成分を$r$個毎, あるいは$l$ を$t$ 毎に区切って考える. 
%区切った$r$ 個の成分からなる分割に対して, その共役分割を考えれば, 
%（最後の区は成分を$r$ 個持たない場合があり, それについて後ほど. ）
%これは$\mathcal{RP}$の定義から成分が$s_1, s_2, \ldots , s_{t-1}, s_t=r$ からなる分割となっており, かつ, それらをどれも少なくとも１つは含む. 
Place the highest number, excluding $r$, on the far right, 
%Place one $s_{t-1}$ on the far right. 
And arrange the rest parts in ascending order. 
We construct tableau by folding and arranging these sequences. \\
%$s_{t-1}$ を一つ除き, 他は順番に並べておく. 
%これをフック状に折り曲げたものを並べて, $\mathcal{CP}$ の分割と対応する盤を回復する. \\
{\textbf{Example7. }}
For $\mathcal{RP}_{\{6;1,3\}}\ni\lambda=(12^110^29^37^16^23^32^11^2)$, 
Let separate parts by $6$, $\lambda^1=(12^110^29^3)$, $\lambda^2=(7^16^23^3), \lambda^3=(2^11^2)$. \\
For $\lambda^1=(12^110^29^3)$, ${}^t({\lambda^1})=(6^93^11^2)$. 
We arrange this parts $116666666663$ 
(Let place the highest number, excluding 6, on the far right. ).
%(let one $s_{t-1}=3$ on the far right. ). 
% $116666666663$ と並べておく. (基本昇順で, $s_{t-1}=3$ を一つだけ右端に. )
Simillary from ${}^t({\lambda^2})=(6^33^31^1), {}^t({\lambda^3})=(3^11^1)$, we get two sequence $1336663, 13$. 
%同様に, ${}^t{\lambda^2}=(6^33^31^1)$から$1336663$, ${}^t{\lambda^3}=(3^11^1)$から$13$ という列が得られる. 
We fold these sequence and construct tableau. 
The question is where to fold. 
It will be uniquely determined if the bottom 6 comes where the number on the right is.
If you fold it in any other way, a row with only 6 or a row with multiple 1s or 3s will appear. \\
%これらをフックに折り曲げ, ヤング盤の形に並べれば良い. 

%その際どこで折り曲げるかが問題となるが, 一番下の$r=6$ が, 丁度右側の$1$や$3$ の行になる様に折り曲げればユニークに決まる. 図でいうと, 仕切りの左が$6$から始まるようにすればよい. 
%（それ以外で折り曲げてしまうと, 行を読んで分割に直す際に, $r=6$ の倍数が出てきてしまったり, 1や3が複数含まれる行が出てきてしまったりする. ）\\
%WinTpicVersion2.15
\unitlength 0.1in
\begin{picture}(48.00,20.00)(4.00,-24.00)
\put(12.5000,-5.5000){\makebox(0,0)[lb]{$6$}}%
\put(14.5000,-5.5000){\makebox(0,0)[lb]{$6$}}%
\put(16.5000,-5.5000){\makebox(0,0)[lb]{$6$}}%
\put(18.5000,-5.5000){\makebox(0,0)[lb]{$3$}}%
\put(12.5000,-7.5000){\makebox(0,0)[lb]{$6$}}%
\put(14.5000,-7.5000){\makebox(0,0)[lb]{$6$}}%
\put(16.5000,-7.5000){\makebox(0,0)[lb]{$3$}}%
\put(12.5000,-9.5000){\makebox(0,0)[lb]{$6$}}%
\put(14.5000,-9.5000){\makebox(0,0)[lb]{$6$}}%
\put(16.5000,-9.5000){\makebox(0,0)[lb]{$3$}}%
\put(12.5000,-11.5000){\makebox(0,0)[lb]{$6$}}%
\put(14.5000,-11.5000){\makebox(0,0)[lb]{$6$}}%
\put(16.5000,-11.5000){\makebox(0,0)[lb]{$1$}}%
\put(12.5000,-13.5000){\makebox(0,0)[lb]{$6$}}%
\put(14.5000,-13.5000){\makebox(0,0)[lb]{$3$}}%
\put(12.5000,-15.5000){\makebox(0,0)[lb]{$6$}}%
\put(14.5000,-15.5000){\makebox(0,0)[lb]{$3$}}%
\put(12.5000,-17.5000){\makebox(0,0)[lb]{$6$}}%
\put(14.5000,-17.5000){\makebox(0,0)[lb]{$1$}}%
\put(12.5000,-19.5000){\makebox(0,0)[lb]{$1$}}%
\put(12.5000,-21.5000){\makebox(0,0)[lb]{$1$}}%
\special{pn 8}%
\special{pa 1200 2200}%
\special{pa 1200 400}%
\special{fp}%
\special{pa 1200 400}%
\special{pa 2000 400}%
\special{fp}%
\special{pa 1800 600}%
\special{pa 1400 600}%
\special{fp}%
\special{pa 1400 600}%
\special{pa 1400 1800}%
\special{fp}%
\special{pa 1600 1200}%
\special{pa 1600 800}%
\special{fp}%
\special{pa 1600 800}%
\special{pa 1800 800}%
\special{fp}%
% STR 2 0 3 0
% 3 400 1700 400 1800 2 0
% good
\put(4.0000,-14.0000){\makebox(0,0)[lb]{good}}%
% STR 2 0 3 0
% 3 2400 1700 2400 1800 2 0
% bad
\put(24.0000,-13.5000){\makebox(0,0)[lb]{bad}}%
\put(30.5000,-5.5000){\makebox(0,0)[lb]{$6$}}%
\put(32.5000,-5.5000){\makebox(0,0)[lb]{$6$}}%
\put(30.5000,-7.5000){\makebox(0,0)[lb]{$6$}}%
\put(32.5000,-7.5000){\makebox(0,0)[lb]{$6$}}%
\put(34.5000,-7.5000){\makebox(0,0)[lb]{$3$}}%
\put(30.5000,-9.5000){\makebox(0,0)[lb]{$6$}}%
\put(32.5000,-9.5000){\makebox(0,0)[lb]{$6$}}%
\put(34.5000,-9.5000){\makebox(0,0)[lb]{$3$}}%
\put(30.5000,-11.5000){\makebox(0,0)[lb]{$6$}}%
\put(32.5000,-11.5000){\makebox(0,0)[lb]{$6$}}%
\put(34.5000,-11.5000){\makebox(0,0)[lb]{$1$}}%
\put(30.5000,-13.5000){\makebox(0,0)[lb]{$6$}}%
\put(32.5000,-13.5000){\makebox(0,0)[lb]{$3$}}%
\put(30.5000,-15.5000){\makebox(0,0)[lb]{$6$}}%
\put(32.5000,-15.5000){\makebox(0,0)[lb]{$3$}}%
\put(30.5000,-17.5000){\makebox(0,0)[lb]{$6$}}%
\put(32.5000,-17.5000){\makebox(0,0)[lb]{$1$}}%
\put(30.5000,-21.5000){\makebox(0,0)[lb]{$1$}}%
\put(34.5000,-5.5000){\makebox(0,0)[lb]{$3$}}%
\put(30.5000,-19.5000){\makebox(0,0)[lb]{$6$}}%
\put(30.5000,-23.5000){\makebox(0,0)[lb]{$1$}}
\special{pn 8}%
\special{pa 3000 400}%
\special{pa 3600 400}%
\special{fp}%
\special{pa 3000 400}%
\special{pa 3000 2400}%
\special{fp}%
\special{pa 3600 600}%
\special{pa 3200 600}%
\special{fp}%
\special{pa 3200 600}%
\special{pa 3200 1800}%
\special{fp}%
\special{pa 3600 800}%
\special{pa 3400 800}%
\special{fp}%
\special{pa 3400 800}%
\special{pa 3400 1200}%
\special{fp}%
\put(42.5000,-5.5000){\makebox(0,0)[lb]{$6$}}%
\put(44.5000,-5.5000){\makebox(0,0)[lb]{$6$}}%
\put(46.5000,-5.5000){\makebox(0,0)[lb]{$6$}}%
\put(48.5000,-5.5000){\makebox(0,0)[lb]{$6$}}%
\put(50.5000,-5.5000){\makebox(0,0)[lb]{$3$}}%
\put(42.5000,-7.5000){\makebox(0,0)[lb]{$6$}}%
\put(42.5000,-9.5000){\makebox(0,0)[lb]{$6$}}%
\put(42.5000,-11.5000){\makebox(0,0)[lb]{$6$}}%
\put(42.5000,-13.5000){\makebox(0,0)[lb]{$6$}}%
\put(42.5000,-15.5000){\makebox(0,0)[lb]{$6$}}%
\put(42.5000,-17.5000){\makebox(0,0)[lb]{$1$}}%
\put(42.5000,-19.5000){\makebox(0,0)[lb]{$1$}}%
\put(44.5000,-7.5000){\makebox(0,0)[lb]{$6$}}%
\put(46.5000,-7.5000){\makebox(0,0)[lb]{$3$}}%
\put(44.5000,-9.5000){\makebox(0,0)[lb]{$6$}}%
\put(44.5000,-11.5000){\makebox(0,0)[lb]{$6$}}%
\put(44.5000,-13.5000){\makebox(0,0)[lb]{$3$}}%
\put(44.5000,-15.5000){\makebox(0,0)[lb]{$3$}}%
\put(44.5000,-17.5000){\makebox(0,0)[lb]{$1$}}%
\put(46.5000,-9.5000){\makebox(0,0)[lb]{$3$}}%
\put(46.5000,-11.5000){\makebox(0,0)[lb]{$1$}}%

\special{pn 8}%
\special{pa 5200 400}%
\special{pa 4200 400}%
\special{fp}%
\special{pa 4200 400}%
\special{pa 4200 2000}%
\special{fp}%
\special{pa 4400 1800}%
\special{pa 4400 600}%
\special{fp}%
\special{pa 4400 600}%
\special{pa 4800 600}%
\special{fp}%
\special{pa 4800 800}%
\special{pa 4600 800}%
\special{fp}%
\special{pa 4600 800}%
\special{pa 4600 1200}%
\special{fp}%
\end{picture}%
\\
Then, $S^{-1}(12^110^29^37^16^23^32^11^2)=(21^115^213^19^27^11^2)\in {\mathcal{CP}}_{\{6;1,3\}}$. 
It's just a coincidence that the exponent is repeats $12$. 
%指数が12でループしているのは偶然. 
\end{proof}
Last, although we have not reached our goal, we will introduce the extension of the theorem to the extent that the correspondence can be proved by the Bessenrodt's map and the inverse map.
%そこには届いていないが, Bssenrodt の対応によって何とかなる部分いっぱいまでの拡張を紹介する
\[
\mathcal{BCP}_{r}(n):=
\left\{
\lambda\in \mathcal{CP}_r(n)\ {\Big{|}}\ 
\begin{array}{l}
{}^\forall i, {\textrm{s.t.}}\ \lambda_i>r(i-1), 
{}^\forall k, \\
Q_r(\lambda_k)=i-1\Rightarrow R_r(\lambda_k)\leq R_r(\lambda_i)
\end{array}
\right\}
\]
\begin{eqnarray*}
&&\mathcal{BRP}_r(n)\\
&:=&
\left\{
\lambda\in \mathcal{RP}_r(n)\ {\Big{|}}\ 
\begin{array}{l}
{}^\forall i, \lambda_{ri}>\lambda_{ri+1},\\
\lambda_{ri}=\lambda_{ri+1}+1\Rightarrow
m_{\lambda_{ri}}\leq \max\{m_{\lambda_{r(i+1)}}, r(i+1)-\ell(\lambda)\}
\end{array}
\right\}
\end{eqnarray*}
\begin{thm}For all $r, n$, 
\[
\sharp\mathcal{BCP}_{r}(n)
=\sharp\mathcal{BRP}_{r}(n)
\]
\end{thm}
Our goal is the correspondence between $\mathcal{CP}_{r}(n), \mathcal{RP}_{r}(n)$ that include the periodic case introduced in this paper.
For that purpose, it is necessary to add some rule to the map used in the proof, and to correspond the partitions not included in ${\mathcal{BCP}}_r$ and ${\mathcal{BRP}}_r$.
%$\mathcal{CP}_{r}(n), \mathcal{RP}_{r}(n)$ 間の対応であり, 本論文で紹介した周期的な場合の対応を精密化として含むものを目標としている. 
%そのためには証明に用いた写像に何かしらのルールを加え, ${\mathcal{BCP}}_r, {\mathcal{BRP}}_r$ に含まれない所を対応させる必要がある. 


\begin{thebibliography}{7}
\bibitem{B}
C. Bessenrodt,
A bijection for Lebesgue's partition identity in the spirit of Sylvester, 
{\it Discrete Math.} {\bf{132}} (1994) 1-10.

\bibitem{S}
J.J. Sylvester, A constructive theory of partitions, arranged in three acts, an interact and an exodion, 
{\it Amer. J. Math.}{\bf{5}}(1882) 251-330.
\end{thebibliography}
\end{document}